              \def\version{12 November, 2010}                    %

\documentclass[reqno,11pt]{amsart}

\usepackage{amsmath, latexsym, srcltx, amssymb}

\newfam\Bbbfam
\font\tenBbb=msbm10
\font\sevenBbb=msbm7
\font\fiveBbb=msbm5
\textfont\Bbbfam=\tenBbb
\scriptfont\Bbbfam=\sevenBbb
\scriptscriptfont\Bbbfam=\fiveBbb

\newcommand{\R}     {\mathbb{R}}
\newcommand{\Z}     {\mathbb{Z}}
\newcommand{\N}     {\mathbb{N}}
\renewcommand{\P}   {\mathbb{P}}

\newcommand{\E}     {\mathbb{E}}
\newcommand{\ssup}[1] {{\scriptscriptstyle{({#1}})}}

\def\1{{\mathchoice {1\mskip-4mu\mathrm l}      
{1\mskip-4mu\mathrm l}
{1\mskip-4.5mu\mathrm l} {1\mskip-5mu\mathrm l}}}

\def\comment#1{}
\newtheoremstyle{thm}{2ex}{2ex}{\itshape\rmfamily}{}
{\bfseries\rmfamily}{}{1.7ex}{}

\newtheoremstyle{rem}{1.3ex}{1.3ex}{\rmfamily}{}
{\itshape\rmfamily}{}{1.5ex}{}

\newtheorem{theorem}{Theorem}[section]
\newtheorem{lemma}[theorem]{Lemma}

\newtheorem{step}{STEP}
\theoremstyle{definition}

\newcommand{\en}       {\end{equation}}
\newcommand{\eq}       {\begin{equation}}
\newcommand{\eqry}   {\begin{eqnarray}}
\newcommand{\enqry}   {\end{eqnarray}}
\newcommand{\eqarray}   {\begin{eqnarray}}
\newcommand{\enarray}   {\end{eqnarray}}
\newcommand{\eqarraystar} {\begin{eqnarray*}}
\newcommand{\enarraystar} {\end{eqnarray*}}
\newcommand{\bel}{\begin{lemma}}
\newcommand{\el}{\end{lemma}}
\newcommand{\bes}{\begin{step}}
\newcommand{\es}{\end{step}}
\newcommand{\bea}{\begin{array}}
\newcommand{\ea}{\end{array}}
\newcommand{\bpr}{\begin{proof}}
\newcommand{\epr}{\end{proof}}

\renewcommand{\section}{\secdef\sct\sect}
\newcommand{\sct}[2][default]{\refstepcounter{section}
\vspace{0.5cm}
\setcounter{equation}{0}
\centerline{ 
\scshape \arabic{section}.\ #1}
\vspace{0.3cm}}
\newcommand{\sect}[1]{
\vspace{0.5cm}
\centerline{\large\scshape #1}
\vspace{0.3cm}}

\renewcommand{\subsection}{\secdef \subsct\sbsect}
\newcommand{\subsct}[2][default]{\refstepcounter{subsection}
\nopagebreak
\vspace{0.5\baselineskip}
{\flushleft\bf \arabic{section}.\arabic{subsection}~\bf #1  }
\nopagebreak}
\newcommand{\sbsect}[1]{\vspace{0.1cm}\noindent
{\bf #1}\vspace{0.1cm}}
\newcommand{\heap}[2]{\genfrac{}{}{0pt}{}{#1}{#2}}

\renewcommand{\subsubsection}{%
\secdef \subsubsect\sbsbsect}
\newcommand{\subsubsect}[2][default]{%
\refstepcounter{subsubsection}
\nopagebreak
\vspace{0.1\baselineskip}
\nopagebreak
{\flushleft
\sffamily\slshape
\arabic{section}.\arabic{subsection}.\arabic{subsubsection}
\ %
\sffamily #1\/.}\ }
\newcommand{\sbsbsect}[1]{\vspace{0.1cm}\noindent
{\bf #1}\ }

\renewcommand{\d}{{\rm d}}
\newcommand{\e}{{\rm e}}

\newcommand{\eps}{\varepsilon}

\newcommand{\smfrac}[2]{\textstyle{\frac{#1}{#2}}}

\newcommand{\Acal}   {{\mathcal A }}

\newcommand{\Ical}   {{\mathcal I }}

\newcommand{\Mcal}   {{\mathcal M }}

\begin{document}

\title[Upper tails of self-intersection local times]{\large 
Upper tails of self-intersection local times\\ \medskip
of random walks: survey of proof techniques}
\author[Wolfgang K\"onig]{}
\maketitle
\thispagestyle{empty}
\vspace{-0.5cm}

\centerline{\sc By Wolfgang K\"onig\footnote{Technical University Berlin, Str. des 17. Juni 136,
10623 Berlin, and Weierstrass Institute for Applied Analysis and Stochastics,
Mohrenstr. 39, 10117 Berlin, Germany, {\tt koenig@wias-berlin.de}}}
\vspace{0.4cm}
\centerline{\textit{WIAS Berlin and TU Berlin}}
\vspace{0.8cm}
\centerline{\textit{\version}}
\vspace{0.8cm}

\begin{quote}{\small {\bf Abstract:} The asymptotics of the probability that the self-intersection local time of a random walk  on $\Z^d$ exceeds its expectation by a large amount is a fascinating subject because of its relation to some models from Statistical Mechanics, to large-deviation theory and variational analysis and because of the variety of the effects that can be observed. However, the proof of the upper bound is notoriously difficult and requires various sophisticated techniques. We survey some heuristics and some recently elaborated techniques and results. This is an extended summary of a talk held on the CIRM-conference on {\it Excess self-intersection local times, and related topics} in Luminy, 6-10 Dec., 2010.}

\end{quote}

\vspace{0.8cm}




\setcounter{section}{0}
\section{Introduction and heuristics}\label{Intro}

\noindent We discuss the logarithmic asymptotics for the upper tails of self-intersection local times of random walks on $\Z^d$. This topic has been studied a lot in the last decade, since it is a natural question, and a rich phenemonology of critical behaviours of the random walk arises, depending on the dimension, the intersection parameter, the scale, and the type of the random process. Furthermore, the question is technically difficult to handle, due to bad continuity and boundedness properties of the self-intersection local time. A couple of different techniques for studying self-intersections have been introduced yet, wich turned out to be more or less fruitful in various situations. It is the goal of this note to survey and compare some of the most fruitful techniques used in recent years.

\subsection{Self-intersection local time}\label{sec-SILT}

\noindent Let $(S_n)_{n\in\N_0}$ be a discrete-time simple random walk in $\Z^d$ started from the origin. We denote by $\P$ the underlying probability measure and by $\E$ the corresponding expectation. The main object of this paper is the {\em self-intersection local time} of the random walk. 
In order to introduce this object, we need the local times of the 
random walk at time $n\in\N$,
\begin{equation}\label{loctim}
\ell_n(z)=\sum_{i=0}^{n}\1_{\{S_i=z\}},\qquad \mbox{ for } z\in\Z^d.
\end{equation}
Fix $p\in(1,\infty)$ and consider the $p$-norm of the local times:
\begin{equation}\label{Lambdadef}
\|\ell_n\|_p
=\Big(\sum_{z\in\Z^d}\ell_n(z)^p\Big)^{1/p},\qquad\mbox { for } n\in\N.
\end{equation}
If $p$ is an integer, then, clearly
\begin{equation}\label{selfinter}
\|\ell_n\|_p^p=\sum_{i_1,\dots,i_p=0}^n\1_{\{S_{i_1}=\dots=S_{i_p}\}}
\end{equation}
is equal to the $p$-fold self-intersection local time 
of the walk, i.e., the number of $p$-fold self-intersections. For $p=2$, this is usually called the {\it self-intersection local time}. For $p=1$, $\|\ell_n\|_p^p$ is 
just the number $n+1$, and for $p=0$, it is equal to $\#\{S_0,\dots,S_n\}$, the 
{\it range\/} of the walk. It is certainly also of interest to study $\|\ell_n\|_p^p$ for non-integer
values of $p>1$, see for example~\cite{HKM04}, where this received technical importance.  The typical behaviour of $\|\ell_n\|_p^p$ has been identified as
\begin{equation}\label{typical}
\E[\|\ell_n\|_p^p]\sim  C a_{d,p}(n),\qquad \mbox{where}\qquad a_{d,p}(n)=
\begin{cases}
n^{(p+1)/2}&\mbox{if }d=1,\\
n (\log n)^{p-1}&\mbox{if }d=2,\\
n&\mbox{if }d\geq 3,
\end{cases}
\end{equation}
for some $C=C_{d,p}\in(0,\infty)$; see \cite{C07} for $d=2$ and \cite{BK07} for $d\geq 3$. In the following, we will concentrate on $d\geq 2$.

\subsection{The problem}\label{sec-problem}

\noindent We are interested in the logarithmic asymptotics of 
$$
\P(\|\smfrac1n \ell_n\|_p\geq r_n),\qquad n\to\infty,
$$
for scale functions $(r_n)_{n\in\N}$ satisfying $r_n-\E[\|\frac 1n\ell_n\|_p]\to\infty$. In order to avoid trivialities and because $\|\frac 1n\ell_n\|_p\leq 1$, we also assume that $r_n<1$. If even $(nr_n)^p\gg a_{d,p}(n)$ (we write $\gg$ if the quotient diverges), we speak of {\it very large deviations}, and if $(nr_n)^p\sim\gamma a_{d,p}(n)$ with $\gamma>1$, we speak of {\it large deviations}. In this note, we will be mainly interested in very large deviations.

In other words, we would like to understand how likely it is for the path to produce many self-intersections, and, additionally, what the typical behavior of the path is on the event $\{\|\frac1n \ell_n\|_p\geq r_n\}$. Certainly, the answer will depend strongly on various issues, like the dimension, the decay of $r_n$, the value of $p$ etc. There is a competition between two effects: clumping together on a small region and the spread-out strength coming from the diffusion mechanism. In order to find the answer, we have to quantify the probabilistic cost of the clumping.

\subsection{Rough heuristics}\label{sec-roughheur}

\noindent Let us give a rough heuristic about what to expect. We consider the very-large deviation case $(n r_n)^p\gg a_{d,p}(n)$. 

The starting point of our heuristic is that the optimal strategy of the path to meet the event $\{\|\frac 1n \ell_n\|_p\geq r_n\}$ is that the path fills a ball $B_{\alpha_n}$ of radius $1\ll\alpha_n\ll n^{1/d}$ within a time interval $[0,t_n]$ with $1\ll t_n\leq n$ in order to produce the required amount $(nr_n)^p$ of self-intersections, and that the path runs freely in the time interval $[t_n,n]$, where he produces the ordinary amount of self-intersections, which is negligible with respect to $(nr_n)^p$. Certainly, the short-time clumping may also take place at some other time instant, e.g.~in the interval $[n-t_n,n]$ or can be divided into several time stretches, but this should not affect the logarithmic asymptotics. We may assume that all the local times $\ell_n(z)\approx \ell_{t_n}(z)$ with $z\in B_{\alpha_n}$ are of the same order, as non-homogeneous strategies are more costly. This order must be equal to $t_n \alpha_n^{-d}$ since altogether $t_n$ hits are distributed on $\#B_{\alpha_n}$ sites. Furthermore, the $p$-norm of the local times is required to equal $nr_n$, i.e.,
$$
(nr_n)^p\asymp \|\ell_n\|_p^p\approx\sum_{z\in B_{\alpha_n}}\ell_{t_n}(z)^p\asymp \alpha_n^d (t_n \alpha_n^{-d})^p=t_n^p \alpha_n^{d(1-p)}, \qquad\mbox{i.e.,}\qquad t_n\asymp nr_n\alpha_n^{d(p-1)/p}.
$$ 
This requires that $\alpha_n\leq r_n^{\frac p{d(1-p)}}$. The negative logarithm of the probabilistic cost to squeeze a $t_n$-step random walk into a ball with radius $\alpha_n$ is of order 
\begin{equation}\label{neglogprob}
-\log\P(S_{[0,t_n]}\subset B_{\alpha_n})\asymp\frac{t_n} {\alpha_n^{2}}\asymp n r_n \alpha_n^{\frac dp(p-1)-2},
\end{equation}
as may be seen from a decomposition of the path into $t_n \alpha_n^{-2}$ equally long pieces of length $\alpha_n^2$ and invoking the central limit theorem. This assertion holds as long as the quantity in \eqref{neglogprob} diverges.

Recall that we want to argue which choice of $\alpha_n$ (i.e., of $t_n$) yields the maximal probability, i.e., the minimal value in \eqref{neglogprob}. It is obvious that this depends on the dimension. Indeed, in the subcritical dimensions $d<\frac {2p}{p-1}$, the exponent of $\alpha_n$ in the last term is negative and hence the optimal choice is to pick $\alpha_n$ (and hence $t_n$) as large as possible. Accordingly, in the supercritical dimensions $d>\frac {2p}{p-1}$, they must be picked as small as possible. Because of the restrictions $t_n\leq n$ and $\alpha_n\geq 1$, this means that the optimal choices are
\begin{equation}\label{tnalphan}
t_n\asymp \begin{cases}
n&\mbox{if }d<\frac {2p}{p-1},\\
n r_n&\mbox{if  }d>\frac {2p} {p-1},
\end{cases},\qquad\mbox{and}\qquad \alpha_n\asymp
\begin{cases}
r_n^{\frac p{d(1-p)}}&\mbox{if }d<\frac {2p}{p-1},\\
1&\mbox{if  }d>\frac {2p} {p-1}.
\end{cases}
\end{equation}
Hence, we may expect that
\begin{equation}\label{roughanswer}
-\frac1n\log\P(\|{\smfrac 1n}\ell_n\|_p\geq r_n)\asymp\frac 1n\frac{t_n} {\alpha_n^{2}}\asymp \frac 1 {\alpha_n^{2}}\asymp
\begin{cases}
r_n^{\frac {2p}{d(p-1)}}&\mbox{if }d<\frac {2p}{p-1},\\
 r_n&\mbox{if  }d>\frac {2p} {p-1}.
\end{cases}
\end{equation}
Recall that $(nr_n)\gg a_{d,p}(n)$, which implies that $r_n\gg n^{\frac 1p-1}$. Hence, $\alpha_n\ll t_n^{1/d}$. This means that the random walk should stay within a region with volume $\ll t_n$ until time $t_n$, and each local time in that region should be of order $t_n/\alpha_n^d\gg 1$.

This heuristics says nothing about the critical dimension $d=\frac {2p}{p-1}$, and this question is indeed delicate to answer (see Section~\ref{sec-Dynkin}). Interestingly, for $r_n\ll 1$, the two scales on the right of \eqref{roughanswer} differ here, but they coincide for $r_n\asymp 1$, where the probability under interest runs on the exponential scale $n$ in any dimension.

\subsection{Precise heuristics}\label{sec-precheur}

\noindent  We now give a heuristic for a more precise version of \eqref{roughanswer}, which strengthens \lq$\asymp$\rq\ to \lq$\sim$\rq\ with explicit identification of the prefactor. This is based on Donsker-Varadhan large-deviation theory. We keep the assumption $ a_{d,p}(n)\ll (n r_n)^p$ (the very-large deviation case) and assume that $r_n\ll 1$ and first turn to the subcritical dimensions $d<\frac{2p}{p-1}$.

Define the scaled normalized version $L_n\colon\R^d\to[0,\infty)$ of the local times $\ell_n$ by
\begin{equation}\label{Lndef}
L_n(x)=\frac{\alpha_n^d}n \ell_n\bigl(\lfloor x\alpha_n\rfloor\bigr),
\qquad\mbox { for }  x\in \R^d.
\end{equation}
Then $L_n$ is a random element of the set of all probability densities on $\R^d$. In the spirit of the celebrated large-deviation
theorem of Donsker and Varadhan, if $\alpha_n$ satisfies
$1\ll \alpha_n^d\ll a_{d,0}(n)$ (see \eqref{typical}), 
then $L_n$ satisfies a weak large-deviation principle 
in the weak $L^1$-topology with speed $n\alpha_n^{-2}$ and rate 
function $\Ical \colon L^2(\R^d)\to[0,\infty]$ given by
\begin{equation}\label{Idef}
\Ical(f)=\begin{cases} \frac 1{2} \bigl\Vert\nabla{f} \bigr\Vert_2^2
&\mbox{if } {f}\in H^1(\R^d)\mbox{ and }\|f\|_2=1,\\ 
\infty&\mbox{otherwise.} 
\end{cases}
\end{equation}
Roughly, this large-deviation principle says that,
\begin{equation}\label{LDPappr}
\P(L_n\in \,\cdot\,)=
\exp\Bigl\{-\frac n{\alpha_n^2} \Bigl[\inf_{f^2\in\,\cdot}\Ical(f)+o(1)\Bigr]\Bigr\},
\end{equation} 
and the convergence takes place in the weak topology. This principle has been partially proved in a special case in \cite{DV79}, a proof in the general case was given in \cite[Prop.~3.4]{HKM04}.

Now note that
$$
\|\ell_n\|_p=\Big(\sum_{z\in\Z^d}\ell_n(z)^p\Big)^{1/p}
= n\alpha_n^{-d}\Big(\sum_{z\in\Z^d}L_n\bigl({\textstyle{\frac z{\alpha_n}}}\bigr)^p\Big)^{1/p}
= n\alpha_n^{\frac{d(1-p)}{p}}\|L_n\|_p=n r_n \|L_n\|_p.
$$
By our choice of $\alpha_n$ in \eqref{tnalphan} with \lq$\asymp$\rq\ replaced by \lq$=$\rq, we have that
\begin{equation}\label{Lambdaident}
\{\|{\smfrac 1n}\ell_n\|_p\geq r_n\}=\big\{\|L_n\|_p\geq 1\big\}\qquad\mbox{and}\qquad \frac n{\alpha_n^2}=n r_n^{\frac {2p}{d(p-1)}}.
\end{equation}
Using \eqref{LDPappr} for the set $\{f\colon \|f^2\|_p\geq 1\}$, we see that
\begin{equation}\label{preciseasymp}
\lim_{n\to\infty}\frac{r_n^{\frac {2p}{d(1-p)}}}n\log\P(\|{\smfrac 1n}\ell_n\|_p\geq r_n)=-\chi_{d,p},
\end{equation}
where 
\begin{equation}\label{chisub}
\chi_{d,p}= \inf\Big\{\frac12 \|\nabla f\|_2^2\colon f\in H^1(\R^d),\|f^2\|_p=1=\|f\|_2\Big\}.
\end{equation}
It turned out in \cite[Lemma~2.1]{GKS04} that $\chi_{d,p}$ is positive if and only 
if $d(p-1)\leq 2p$. Formula \eqref{preciseasymp} is the precise version of \eqref{roughanswer}. We see that the main contribution to the event $\{\|\frac 1n\ell_n\|_p\geq r_n\}$ comes from those random walk realisations that make the rescaled local times, $L_n$, look like the minimiser(s) $f^2$ of the variational formula on the right-hand side of \eqref{chisub}.

An analogous heuristic applies for the supercritical dimensions $d>\frac {2p}{p-1}$. We keep the assumption $ a_{d,p}(n)\ll (n r_n)^p$, but drop the assumption that $r_n\ll 1$ and assume that $r=\lim_{n\to\infty}r_n\in[0,1]$ exists. Pick $\alpha_n=1$, and the time $n$ must be reduced to $s t_n=snr_n\leq n$, and later we optimize over $s\in(0,1/r)$. Hence, we approximate 
$$
\{\|{\smfrac 1n}\ell_n\|_p\geq r_n\}\approx \big\{\|\ell_{s t_n}\|_p\geq nr_n\big\}=\Big\{\Big\|\frac 1{s t_n}\ell_{s t_n}\Big \|_p \geq \frac 1s\Big\}.
$$
(The set is non-empty only for $s\geq 1$, but this will come out naturally when optimising.) This time we use that $\frac 1{s t_n}\ell_{s t_n}$ satisfies a large-deviation principle on the set of probability measures on $\Z^d$ with scale $st_n$. The rate function $\Ical^{\ssup{\rm d}}$ comes via a contraction principle from a principle for the empirical measures of Markov chains; we abstain from writing it down. Hence, we obtain
\begin{equation}\label{preciseheursuper}
\lim_{n\to\infty}\frac1{nr_n}\log \P(\|{\smfrac 1n}\ell_n\|_p\geq r_n)=-\chi_{d,p},
\end{equation}
where
\begin{equation}\label{chisuper}
\begin{aligned}
\chi_{d,p}&=\inf_{s\in(0,1/r)}s\inf\Big\{\Ical^{\ssup{\rm d}}(g^2)\colon g\in\ell^2(\Z^d), \|g^2\|_p= \frac 1s, \|g\|_2=1\Big\}\\
&=\inf\Big\{\frac{\Ical^{\ssup{\rm d}}(g^2)}{\|g^2\|_p}\colon g\in\ell^2(\Z^d),\|g\|_2=1,\|g^2\|_p\geq r\Big\}.
\end{aligned}
\end{equation}
We see that the main contribution to the event $\{\|\frac 1n\ell_n\|_p\geq r_n\}$ comes from those random walk realisations that make the normalized local times $\frac 1{s t_n}\ell_{s t_n}$ equal to a minimizer $g^2$ of the formula on the right-hand side of \eqref{chisuper} inside some box of bounded radius. In particular, the parameter $s$ should be chosen as $\|g^2\|_p^{-1}$. After time $st_n$, the random walk leaves that bounded box and runs like a free simple random walk and produces a negligible amount of self-intersections.

\subsection{Continuous-time random walks}\label{sec-Conttime}

\noindent The assertion in \eqref{preciseasymp}-\eqref{chisub} should also be literally true for a continuous-time simple random walk $(S_t)_{t\in[0,\infty)}$, and also the scale in \eqref{preciseheursuper} should be the same as in the discrete-time case. However the rate function (and therefore the formula for $\chi_{d,p}$) is different in the super-critical dimensions: it is  $g^2\mapsto \frac 12\|\nabla g\|_2^2=\frac 12 \sum_{z,z'\in\Z^d\colon z\sim z'}(g(z)-g(z'))^2$, which is the discrete-space version of the principle that $L_n$ satisfies. Using a simple scaling argument we have that \eqref{chisuper} must be replaced here by
\begin{equation}\label{chisupercont}
\begin{aligned}
\chi_{d,p}=\inf\Big\{\frac{\frac 12\|\nabla g\|_2^2}{\|g^2\|_p}\colon g\in\ell^2(\Z^d),\|g\|_2=1,\|g^2\|_p\geq r\Big\}
\end{aligned}
\end{equation}
which reduces in the case $r_n\to r=0$ to
\begin{equation}\label{chisupercontr=0}
\chi_{d,p}=\inf\Big\{\frac 12\|\nabla g\|_2^2\colon \|g^2\|_p=1\Big\}.
\end{equation}

\subsection{Exponential moments}\label{sec-expomom}

\noindent The statements in \eqref{preciseasymp} and \eqref{preciseheursuper} are intimately connected with analogous statements about the exponential moments of $\|\ell_n\|_p$. This is a version of the well-known G\"artner-Ellis theorem (see  \cite[Sect.~4.5]{DZ98}). More precisely, via the exponential Chebyshev inequality, they follow from the logarithmic asymptotics of suitable exponential moments, and the lower bound can be proved with the help of a transformation in the spirit of the Cram\'er transform. 

More precisely, if $d<\frac{2p}{p-1}$, abbreviate  $\lambda=\frac{2p+d-dp}{2p}\in(0,1)$, then \eqref{preciseheursuper} follows from the assertion
\begin{equation}\label{ExpMom1}
\frac{1}{n}\log\E\Big( \e^{\theta_n\| \ell_n\|_p}\Big)\sim
\theta_n^{1/\lambda}\rho_{d,p}^{\ssup{\rm c}}(1),\qquad n\to\infty,
\end{equation}
for $(a_{d,p}(n)^{1/p}/n)^{\lambda/(1-\lambda)}\ll\theta_n\ll1$, where
\begin{equation}\label{rhocdef}
\begin{aligned}
\rho_{p,d}^{\ssup{\rm c}}(\theta) &= \sup\Big\{\theta\|f^2\|_p-\frac{1}{2}\|\nabla f\|_2^2\colon f\in H^1(\R^d), \|f\|_2=1\Big\}\\
&=\theta^{1/\lambda}\lambda\left(\frac{2p}{d(p-1)}\chi_{d,p}\right)^{\frac{\lambda-1}{\lambda}},\qquad\theta>0.
\end{aligned}
\end{equation}
Indeed, apply the exponential Chebyshev inequality with $\theta_n=(r_n\lambda/\rho_{d,p}^{\ssup {\rm c}}(1))^{\lambda/(1-\lambda)}$ and to use the second line of \eqref{rhocdef} (which can be shown elementarily by scaling arguments, see \cite[Remark 1.3]{BK10}), to derive the upper bound in \eqref{preciseasymp} . The reason that also the lower bound can be shown with the help of a Cram\'er-type transformation using \eqref{ExpMom1} is that $-\log \P(\|\frac 1n \ell_n\|_p\geq r_n)$ is asymptotically a convex function of $r_n$ (it is $\sim \chi_{d,p} n r_n^{\frac{2p}{d(p-1)}}$, and the power is larger than one); note that this method, the G\"artner-Ellis method, produces only convex rate functions.

In the supercritical dimension, this line of arguments works as well in the case $r_n\ll 1$ since  $-\log \P(\|\frac 1n \ell_n\|_p\geq r_n)$ is asymptotically linear in $r_n$. However, in the case $r_n\to r\in(0,1)$, it does not seem to work since both $\chi_{d,p}$'s defined in \eqref{chisuper} and in \eqref{chisupercont} depend on $r$, and it seems not clear whether the map $r\mapsto r\chi_{d,p}$ is convex. This is also reflected by the fact that the logarithmic rate of the exponential moments of $\|\ell_n\|_p$ is possibly not differentiable, see \cite[Theorem~1.1(i) and Remark~1.5]{BK10}, where it was shown that, for any $\theta>0$, for continuous-time random walk,
\begin{equation}\label{DiscreteExpMom}
\lim_{t\rightarrow \infty}\frac{1}{t}\log\E\Big( \e^{\theta\| \ell_t\|_p }\Big)=\rho_{p,d}^{\ssup{\rm d}}(\theta)=\sup\left\{\theta\| g^2\|_p-\frac12 \|\nabla g\|_2^2\colon g\in\ell^2(\Z^d),g\geq 0,\|g\|_2=1 \right\}.
\end{equation}
Note that the right-hand side is the discrete version of $\rho_{p,d}^{\ssup{\rm c}}(\theta)$ defined in \eqref{rhocdef}.

\subsection{Difficulties}\label{sec-difficult}

\noindent There are several serious obstacles to be removed when trying to turn the above 
heuristics into an honest proof: (1) the large-deviation principles only hold on 
compact subsets of $\R^d$ resp.~$\Z^d$, (2) the functional $f^2\mapsto \|f^2\|_p$ is not bounded in continuous, nor in discrete space, and (3) this functional is not continuous in the topology of the large-deviation principle. 

Removing the 
obstacle (1) is easy and standard (see Section~\ref{sec-compact}), but it is in general notoriously difficult to 
overcome the obstacles (2) and (3) for related problems. In the subcritical dimensions, the transition from discrete to continuous space while taking the limit causes additional technicalities. In the supercritical dimensions, the reduction of the time scale from $n$ to $s t_n$ is also hard to justify rigorously. The critical dimension $d=\frac{2p}{p-1}$, i.e., $p=\frac d{d-2}$, is even more delicate since the question if the discrete or the continuous picture arises seems to depend on the precise choice of $r_n$. See Section~\ref{sec-Dynkin} for a rigorous answer.

These difficulties make the proofs of \eqref{preciseasymp} and \eqref{preciseheursuper} a demanding task.

\subsection{Compactification}\label{sec-compact}

\noindent In most of the proofs of upper bounds for probabilities under interest here, one of the main steps is to estimate $\|\ell_n\|_p\leq \|\ell_n^{\ssup R}\|_p$, where $\ell_n^{\ssup R}$ are the local times of the periodized version of the walk  in the box $B_R=[-R,R]^d\cap\Z^d$ with $R=R_n=L\alpha_n$ and a large parameter $L$. This estimate is easily verified and understood: when putting the free walk onto the torus, one does not lower the number of self-intersections, but possibly increases them. Hence, one is left with the same task for the periodized walk, which lives on a compact part of the space $\Z^d$, which depends on $n$. If one can manage the problem on the torus $B_{L\alpha_n}$ up to logarithmic equivalence, one ends up with an $L$-dependent variational formula, which is elementarily shown to converge towards the correct one as $L\to\infty$. In these notes, we will therefore sometimes tacitly impose the condition $S_{[0,n]}\subset B_R$ without mentioning that the transition probabilities of the walk have been slightly changed. However, this estimate is useful only in the cases in which the typical behavior of the path is to fill a centred box of side length $R$ more or less uniformly. This applies to the subcritical dimensions, but rules out the supercritical ones.

\subsection{Lower bounds}\label{sec-lowbound}

\noindent Actually, the proof of the lower bounds in \eqref{preciseasymp} and \eqref{preciseheursuper} is quite simple and is done as follows in the subcritical dimensions. Pick $q>1$ such that $\frac 1p+\frac 1q=1$ and pick some continuous and bounded function $f$ having compact support and satisfying $\|h\|_q=1$. Then H\"older's inequality gives that $\|L_n\|_p\geq \langle h,L_n\rangle$. Now the large-deviation principle for $L_n$ can safely be applied to $\langle h,L_n\rangle$, since the map $f^2\mapsto \langle h, f^2\rangle $ is continuous and bounded in the topology of the principle. Hence, we obtain the lower bound in \eqref{preciseasymp} with $\chi_{d,p}$ replaced by $\inf\{ \Ical(f)\colon f^2\in L^p(\R^d),\langle h,f^2\rangle=1\}$. Optimizing over $h$ and thereby using the duality between $L^p$ and $L^q$, we see that the lower bound of \eqref{preciseasymp} arises, after employing some elementary approximation arguments. A similar argument applies in the supercritical dimensions.

\section{Techniques for proving upper bounds}\label{sec-proved}

\noindent In this section, we survey various techniques to prove the upper bound in the statements \eqref{preciseasymp} and \eqref{preciseheursuper} and closely related variants of them.

\subsection{Triangular decomposition and smoothing}\label{sec-Smooting}

\noindent In a long series of papers, among a lot of further results on intersections of random motions, Chen also gives a proof of \eqref{preciseasymp} in the most interesting cases $d=2=p$ and $d=3, p=2$, see \cite[Theorems 8.2.1 and 8.4.2]{Ch09}. Actually, he admits more general random walks and much smaller choices of the scale function $(r_n)_{n\in\N}$. He shows that \eqref{preciseasymp} is even true for $r_n=\frac 1n (\E[\|\ell_n\|_2^2]+n b_n)^{1/2}$ with $1\ll b_n\ll n$. Actually, he proves the exponential version \eqref{ExpMom1}.

The three main ideas of the proof method he uses are a triangular decomposition of the number of self-intersections (this restricts the method to $p=2$), a smoothing technique with the help of a convolution of a smooth approximation of the delta measure, and a series of  Banach space tools including the Minkowski functional, the Hahn-Banach theorem and Arzel\'a-Ascoli's theorem.

Indeed, he writes
\begin{equation}\label{split1}
\|\ell_n\|_2^2=\sum_{j=1}^{2^N} \eta^{\ssup N}_j+\sum_{j=1}^N \sum_{k=1}^{2^{j-1}}\xi^{\ssup N}_{j,k},
\end{equation}
where $N\in\N$ is a large auxiliary parameter and
\begin{equation}\label{split2}
\begin{aligned}
\eta^{\ssup N}_j&=\sum_{(j-1)n2^{-N}<i<i'\leq jn2^{-N}}\1\{S_i=S_{i'}\},\\
\xi^{\ssup N}_{j,k}&=\sum_{\heap{(2k-2)n2^{-j}<i\leq (2k-1)n2^{-j}}{(2k-1)n2^{-j}<i'\leq (2k)n2^{-j}}}\1\{S_i=S_{i'}\}.
\end{aligned}
\end{equation}
This decomposition was already used by Le Gall \cite{Le86}, it can also be defined via an iterated bisection of the path. Its advantage is that $\eta^{\ssup N}_1,\dots,\eta^{\ssup N}_{2^N}$ are i.i.d.~with distribution equal to the number of self-intersections of a random walk of length $\approx n2^{-N}$ and that, for any $j\in\{1,\dots,N\}$, the variables $\xi^{\ssup N}_{j,1},\dots,\xi^{\ssup N}_{j,2^{j-1}}$ are i.i.d.~ with distribution equal to the number of mutual intersections of two independent random walks of length $\approx n 2^{-j}$. This decomposition was already used in the 1960ies for the study of the self-intersections of two-dimensional Brownian motion.

The second idea is to convolute the normalised and rescaled local times $L_n$ defined in \eqref{Lndef} with some smooth approximation, $\varphi_\eps$, of the Dirac delta measure as $\eps\downarrow0$. The replacement of $\|L_n\|_p^p$ with the smoothed ones, $\|L_n\star \varphi_\eps\|_p^p$, with full control of the  asymptotics  requires some technical care, but can be done using more or less standard means. 

The large-deviation arguments for $\|L_n\star \varphi_\eps\|_p$ are easier to derive than for $\|L_n\|_p$, but however require some substantial work, see \cite[Sect.~4.2]{Ch09}. The reason is that the map $L_n\mapsto \|L_n\star \varphi_\eps\|_p$ has still bad continuity properties. Chen's ingenious way to solve this problem uses a compactness criterion introduced in \cite{dA85}, formulated in terms of bounds for certain exponential integrals of the Minkowski functional of a convex, positively balanced set. The way to make this criterion applicable is long and uses a series of ideas from functional analysis, like the Arzel\'a-Ascoli theorem, topological duality between the spaces $L^p$  and $L^q$ for $\frac 1p=\frac 1q=1$, and the Hahn-Banach theorem.

\subsection{Iterated bisection}\label{sec-binsplit}

\noindent 
As we mentioned in Section~\ref{sec-Smooting}, Le Gall \cite{Le86} introduced a technique of successive division of the path into approximately equally long pieces and controlling the self-interaction of each piece and the mutual interaction between them. This induction procedure is equivalent to the splitting technique described in \eqref{split1}-\eqref{split2}. {\it A priori} this method works only for $p=2$. However, it has been further developed by Asselah \cite{A09}  to be used for any value of $p\in(1,\infty)$. This enables him to prove both assertions in \eqref{roughanswer} for both large and very large deviations. However, his approach admits only a study of dimensions $d\geq 3$, since he uses transience of the walk at some place.

The kernel of Asselah's bisection technique for $\|\ell_n\|_p^p$, i.e., for a sum of $p$-th powers of integers, is the estimate
$$
(l_1+l_2)^p\leq l_1^p+l_2^p+2^p\sum_{i=0}^\infty b_{i+1}^{p-2}l_1 l_2\1\{b_i\leq \max\{l_1,l_2\}< b_{i+1}\}\qquad l_1,l_2\in\N,
$$
where $1=b_0<b_1<b_2<\dots$ defines a suitable partitioning of $[1,\infty)$. Using this estimate iteratively for bisections of the path, one obtains an upper bound for $\|\ell_n\|_p^p$ in terms of a sum of the $p$-norms of the respective fragments of the path (which are independent) plus an additional term coming from their mutual interaction. One additional ingredient of the proof is a decomposition of the space into regions where the local times are small, medium-sized or large. The event $\{\|\ell_n\|_p\geq r_n\}$ is decomposed in several partial events, whose probabilities are estimated using various arguments.

\subsection{Surgery on circuits and clusters}\label{sec-surgery}

\noindent As we explained in Section~\ref{sec-precheur}, in the supercritial dimension, the significant contribution to a large value of the intersection local time comes from paths that have extremely high values on a bounded region. This intuitive picture is the leading idea in the proof given in \cite{A08b} (see also \cite{A08a}), where \eqref{preciseheursuper} is proved for $p=2$, $d\geq5$ and $r_n\asymp n^{-1/2}$, i.e., for the large-deviation regime.

The main technical tool is an upper estimate of $\|\ell_n\|_2^2-\E[\|\ell_n\|_2^2]$ in terms of $\|\1_\Lambda \ell_{s\sqrt n}\|_2^2$ for many choices of a finite set $\Lambda\subset\Z^d$ on the event $\{S_{s\sqrt n}=0\}$, i.e., for a circuit. To derive this, Asselah introduces for infinite-time random walk, using some iterative procedure called surgery, a map from finite $n$-dependent boxes to bounded subboxes that compares paths with high values of local times in the large box to those having high local time values in the small box. Particular attention is given to the region where the local times are of order $\sqrt n$; finally it is shown that this set is bounded in $n$.

The outcome of this technique is that the existence and non-triviality of the limit in  \eqref{preciseheursuper} is shown. In a second step, its value is identified as $\frac 12$ times the constant on the right-hand side of \eqref{CM091} by a comparison between the two problems of mutual intersections of two independent walks and self-intersections of one walk.

\subsection{Dynkin's isomorphism}\label{sec-Dynkin}

\noindent The critical choice $p=\frac d{d-2}$ in dimensions $d\geq 3$ is considered in \cite{C10}. Actually, it is shown there that, in the continuous-time case, \eqref{preciseasymp} is true with $\chi_{d,p}$ as in \eqref{chisupercont}, for any $n^{\frac 1p-1}\ll r_n\ll 1$. This interestingly shows that the critical dimension $d=\frac{2p}{p-1}$ belongs to the lower critical case, as it concerns the radius $\alpha_n$ of the ball in which the main bulk of the self-intersections occur, but to the upper critical dimension, as it concerns the nature of the variational formula describing the precise logarithmic asymptotics. 

The main idea used in \cite{C10} is Dynkin's isomorphism theorem \cite{D88a}, which says that the joint law of the local times of a symmetric recurrent Markov process stopped at an independent exponential time is related to the law of the square of a Gaussian process whose covariance function is the Green kernel of the stopped Markov process. To apply this, in a first step, the exponential moments of $\|\ell_t\|_p$ are estimated from above against the exponential moments of $\|\ell^{\ssup R}_\tau\|_p$, where $\ell^{\ssup R}$ are the local times of the torus version of the walk on $B_R$, and $\tau$ is an independent exponential time with parameter $\asymp r_t$. Now introduce a Gaussian process $Z=(Z_x)_{x\in B_R}$ with covariance matrix equal to the Green function, $G_{R,\tau}$, of the stopped walk $(S^{\ssup R}_{t\wedge \tau})_{t\in[0,\infty)}$ on the torus $B_R$. Then the exponential moments of $\|\ell^{\ssup R}_\tau\|_p$ can be written in terms of exponential moments of $\|Z\|_{2p}^2$ with some slightly modified density. The great advantage of this rewrite is that now concentration inequalities for Gaussian integrals can be applied to the exponential moments of $(\|Z\|_{2p}-M)^2$, where $M$ denotes the median of $\|Z\|_{2p}$. These inequalities are so precise that they prove the crucial fact that the tail behaviour of $\|Z\|_{2p}-M$ is equal to that of a Gaussian variable with variance equal to $\sup\{\langle f, G_{R,\tau} f\rangle\colon f\in\ell^{2p}(\Z^d),  \|f\|_{2p}=1\}$. If one picks $R\asymp t^{1/d}$, then this supremum converges towards  $\chi_{d,p}$ defined in \eqref{chisupercont}, and this is the kernel of this proof method.

\subsection{Polynomial moments}\label{sec-polymom}

\noindent Another sucessful technique  is based on an expansion of $\exp\{\theta \alpha_t^{2-d+d/p}\|\ell_t\|_p\}$ and a precise estimation of the polynomial moments of $\|\ell_t\|_p$ with suitable $t$-dependent powers. More precisely, in subcritical dimensions in the proof of  \cite[Prop.~2.1]{HKM04} it is shown that, for any $L\in(0,\infty)$, in the time-continuous case,
\begin{equation}\label{HKM041}
\E\big(\|\ell_t\|_p^{pk}\1\{S_{[0,t]}\subset B_{L\alpha_t}\}\big)\leq k^{kp} C^k \alpha_t^{k[d+(2-d)p]},\qquad k\geq \frac t{\alpha_t^2},
\end{equation}
for some $C\in(0,\infty)$ and for all sufficiently large $t$. It is easy to see that this implies that 
$$
\limsup_{\theta\downarrow 0}\limsup_{t\to \infty}\frac{\alpha_t^2}{t}\log\E\Big[\exp\Big\{\theta\alpha_t^{-\frac1p [d+(2-d)p]}\|\ell_t\|_p\1\{S_{[0,t]}\subset B_{L\alpha_t}\} \Big]\leq 0,
$$
which was one of the partial goals in \cite{HKM04}. This statement is a bit less than the upper bound in \eqref{ExpMom1}, but identifies the correct scale.

The proof of \eqref{HKM041} consists of a combinatorial analysis of the polynomial moments by explicitly writing out $\ell_t(z)=\int_0^t\delta_z(S_r)\,\d r$ and the $pk$-th moments and summarizing and transforming the arising multi-sum as far as possible. No attempt to optimize \eqref{HKM041} nor to find the best value of $C$ was made in \cite{HKM04}. The method works in any subcritical dimension, but only for $\alpha_t\ll t^{1/(d+1)}$, which is a severe restriction. The kernel of the reason why this methods works is the integrability of the $p$-th power of the Green function of Brownian motion around its singularity.

This method was applied also in the supercritical dimensions in \cite{CM09} for the closely related problem of the mutual intersections of $p$ independent copies $S^{\ssup 1},\dots,S^{\ssup p}$ of $(S_n)_{n\in\N_0}$ rather than the self-intersections of one walk. Here it is possible (and the main interest of \cite{CM09}) to consider these intersections with infinite time horizon and to study its upper tails. Denote by 
$$
I=\sum_{i_1,\dots,i_p=0}^\infty\1\{S^{\ssup 1}_{i_1}=\dots=S^{\ssup p}_{i_p}\}
$$ 
this intersection local time, then the main result of \cite{CM09} is
\begin{equation}\label{CM091}
\lim_{a\to\infty}a^{-1/p}\log \P(I>a)=-p\inf\big\{\|h\|_q\colon h\in\ell^q(\Z^d), h\geq 0,\|\Acal_h\|\geq1\big\},
\end{equation}
where $\frac 1p+\frac 1q=1$, and the operator $\Acal_h\colon\ell^2(\Z^d)\to\ell^2(\Z^d)$ is defined by 
$$
\Acal_hg(x)=\sqrt{\e^{h(x)}-1}\sum_{y\in\Z^d}G(x,y)\sqrt{\e^{h(y)}-1},
$$
and $G$ is the Green function of the walk. For continuous-time simple random walk, the right-hand side of \eqref{CM091} is shown to equal $-\chi_{d,p}$ defined in \eqref{chisupercont}. It is expected that for the infinite-time self-intersection local times, i.e., after replacing $I$ by $\|\ell_\infty\|_p^p$, \eqref{CM091} remains true with the factor $p$ removed on the right-hand side, and that $\|\ell_\infty\|_p^p$ may also be replaced by $\|\ell_{f(a)}\|_p^p$ for some explicit function $f(a)$. Similarly, also in the discrete-time case, the right-hand side of of \eqref{CM091} should be equal to $-\chi_{d,p}$ in \eqref{chisuper}. Details are currently being worked out in \cite{BK11+}.

The proof of \eqref{CM091} is again based on the asymptotical analysis of high polynomial moments. It is used \cite[Lemma~2.1]{KM02} that, for any positive random variable $X$,
$$
\lim_{k\to\infty}\frac 1k\log\E\Big[\frac{X^k}{k!^p}\Big]=\kappa\qquad\Longleftrightarrow\qquad\lim_{a\to\infty}a^{-1/p}\log \P(X>a)=-p\e^{\kappa/p}.
$$
For the identification of the high polynomial moments of $I$, some compactification procedure is developed that is in the spirit of the periodization idea mentioned in Section~\ref{sec-compact}, but this time for the $p$-th powers of the Green function of the walk instead of the $p$-th power of the local times. The fact that this procedure gives the correct upper bound may be interpreted by saying that the main bulk of the intersections occur in some box of bounded radius, which may be far from the origin.

\subsection{Density of local times}\label{sec-densityloctim}

\noindent The following is restricted to continuous-time random walk $(S_t)_{t\in[0,\infty)}$. The approach of \cite{BK10} is to employ an explicit formula for the joint density of the local times $(\ell_t(z))_{z\in B}$ in a finite subset $B$ of $\Z^d$, which has been derived in \cite{BHK07}. This makes it possible to explicitly write down a formula for the expected exponential moments of $\|\ell_t\|_p$ on the event $\{S_{[0,t]}\subset B\}$. Even though the representation for this density derived in \cite[Theorem~2.1]{BHK07} is almost impossible to penetrate, \cite[Theorem~3.6]{BHK07} gives a handy upper bound for such expectations.

For the subcritical dimension, in \cite{BK10}, it was obtained in this way that
\begin{equation}\label{BK101}
\frac 1t\log\E\left(\exp\Big\{t\alpha_t^{-2\lambda}\big\|\smfrac 1t\ell_t^{\ssup{L\alpha_t}}\big\|_p\Big\}\right)
\leq  \rho^{\ssup{\rm d}}_{d,p}(L\alpha_t,\alpha_t^{-2\lambda})+\eps_t,
\end{equation}
where we recall that $\lambda=\frac{2p+d-dp}{2p}\in(0,1)$ and 
$$
\rho^{\ssup{\rm d}}_{d,p}(R,\theta)=\sup_{\mu\in \Mcal_1(B_{R})}\left[\theta\|\mu\|_p-\| \left(-A_{R}\right)^{1/2}\sqrt{\mu}\|_2^2 \right],
$$
and $A_{R}$ is the generator of the periodized version of the random walk on the box $B_{R}$, recall Section~\ref{sec-compact}, and $\eps_t$ is some explicit error term. It is required that $\eps_t\leq \exp\{o( r_t)^{\frac {2p}{d(p-1)}}\}$, and this in turn enforces that $r_t\gg({\log t}/t)^{\frac{d(p-1)}{p(d+2)}}$, which imposes a restriction on the validity for technical reasons.

The main term on the right-hand side of \eqref{BK101}, $\rho^{\ssup{\rm d}}_{d,p}(L\alpha_t,\alpha_t^{-2\lambda})$ is a compact-space version of $\rho_{p,d}^{\ssup{\rm d}}(\theta) $ introduced  in \eqref{DiscreteExpMom}, and there is a close connection to the continuous version defined in \eqref{rhocdef}. Indeed, the main work in \cite{BK10} is devoted to the proof of
\begin{equation}\label{BK102}
\limsup_{L\to\infty}\limsup_{t\to\infty}\alpha_t^2\rho^{\ssup{\rm d}}_{d,p}(L\alpha_t,\alpha_t^{-2\lambda})\leq \rho_{p,d}^{\ssup{\rm c}}(1),
\end{equation}
and this finishes the proof of the upper bound in \eqref{ExpMom1}. The proof of \eqref{BK102} is in the spirit of Gamma-convergence techniques, some elements of finite element theory is employed. Unfortunately, in the course of the proof, the technical assumption that $d<\frac2{p-1}$ must be made, which severely restricts the validity in the dimension.

\end{document}